\documentclass[12pt]{article}
\usepackage{amsmath}
\usepackage[noend]{algpseudocode}
\usepackage{algorithm}
\usepackage{amssymb}
\usepackage{graphicx}
\usepackage{psfrag}
\usepackage{hyperref}
\usepackage{enumitem}
\usepackage{subcaption}

\newcommand{\reals}{{\mbox{\bf R}}}

\newcommand{\diag}{\mathop{\bf diag{}}}

\newcommand{\ie}{{\it i.e.}}
\newcommand{\argmin}{\mathop{\rm argmin}}

\newcommand{\ones}{\mathbf 1}

\newcommand{\BIT}{\begin{itemize}}
\newcommand{\EIT}{\end{itemize}}
\newcommand{\BEQ}{\begin{equation}}
\newcommand{\EEQ}{\end{equation}}
\newcommand{\BEAS}{\begin{eqnarray*}}
\newcommand{\EEAS}{\end{eqnarray*}}

\newcounter{algorithmctr}[section]
\renewcommand{\thealgorithmctr}{\thesection.\arabic{algorithmctr}}
{\refstepcounter{algorithmctr}
\begin{list}{}{%
\setlength{\rightmargin}{0\linewidth}%
\setlength{\leftmargin}{.05\linewidth}}%
\rmfamily\small
\item[]{\setlength{\parskip}{0ex}\hrulefill\par%
\nopagebreak{\bfseries\textsf{Algorithm \thealgorithmctr~}}}}%
{{\setlength{\parskip}{-1ex}\nopagebreak\par\hrulefill}
\end{list}}

{\refstepcounter{algorithmctr}
\begin{list}{}{%
\setlength{\rightmargin}{0\linewidth}%
\setlength{\leftmargin}{.05\linewidth}}%
\rmfamily\small
\item[]{\setlength{\parskip}{0ex}\hrulefill\par%
\nopagebreak{\bfseries\textsf{Algorithm}}}}%
{{\setlength{\parskip}{-1ex}\nopagebreak\par\hrulefill}
\end{list}}

\begin{document}

\title{A Distributed Method for Cooperative Transaction Cost Mitigation}
\author{Nikhil Devanathan \and Logan Bell \and Dylan Rueter
\and Stephen Boyd}
\date{\today}
\maketitle

\begin{abstract}
Funds at large portfolio management firms may consist of many portfolio
managers (PMs), each managing a portion of the fund and optimizing a 
distinct objective. Although the PMs determine
their trades independently, the trade lists may be netted and executed by the
firm. These net trades may be sufficiently large to impact the market prices,
so the PMs may realize prices on their trades that are different
from the observed midpoint price of the assets before execution. These
transaction costs generally reduce the returns of a portfolio over time.  We
propose a simple protocol, based on methods from distributed convex
optimization, by which a firm can communicate estimated transaction costs to
its PMs, and the PMs can potentially revise their
trades to realize reduced transaction costs. This protocol does not require the
PMs to disclose their method of determining trades to the firm
or to each other, nor does it require the PMs to communicate
their trade lists with each other.  As the number of adjustment rounds grows,
the trades converge to the ones that are optimal for the firm. As a practical
matter we observe that even just a few rounds of adjustment lead to substantial
savings for the firm and the PMs.
\end{abstract}

\clearpage

\section{Introduction}\label{s-intro}

\paragraph{The coordination problem.}
We consider the setting of a systematic firm with multiple PMs managing
independent investment sleeves (sub-portfolios) of a single fund\footnote{
This paper is a purely mathematical contribution to distributed
optimization in the context of portfolio management. It does not
describe, reference, or prescribe the practices of any specific firm.
By \emph{netting} trade lists, we refer to the mathematical aggregation
of individual trade lists into a single net trade vector. In practice,
the internal crossing of trade lists is a regulated activity; we
exclusively consider the case in which crossing is entirely permissible
and make no assertions about its applicability in any particular
regulatory or operational context.
We assume, as a modeling abstraction, that there exists a central entity
within the firm charged with aggregating the PMs' trade lists and
executing the resulting net trade on the market. This is a simplifying
modeling device.
Finally, we use the term \emph{cooperative} throughout this paper in its
game-theoretic sense, and not
in the colloquial sense of informal coordination or agreement between
parties.}.
Each PM determines their desired trading actions by
solving a specific optimization problem, but the outcomes of these actions
are often coupled at the firm level. This coupling arises when the cost or
utility realized by the firm depends on the aggregate actions of the PMs.

A primary example of this is netted transaction costs. The firm may net the
trade lists of all the PMs and execute the net trade on the market. Trading
impacts market prices, resulting in transaction costs that depend on the
aggregate trades. Similar coupling can arise from shared borrow costs or
firmwide risk constraints. The challenge is to coordinate the PMs to optimize a
global firm objective while respecting the autonomy and private information of
the individual PMs.

\paragraph{Joint transaction cost optimization.}
We assume that all PMs determine their trade lists by solving an optimization
problem, typically convex, that takes into account their views of future risk
and returns, the cost of holding short positions, and the cost of trading, as
well as portfolio specific constraints, such as limits on what assets can be
held. To find the optimal trades taking the netted transaction cost into
account, we can collect all PMs' optimization problems, sum their objectives,
and subtract the transaction cost associated with the net trades. The
individual optimization problems associated with each portfolio are now coupled
through the net trade transaction cost. If the optimization problems used by
each PM are convex, we obtain a large convex optimization problem, which
determines all trades simultaneously, minimizing the total transaction cost to the
firm, in addition to maximizing each PM's objective. We refer to this as the
joint problem, since all the trades are found by solving a single optimization
problem that takes the net trade transaction cost into account.

\paragraph{Distributed solution of the joint problem.}
We propose a distributed method to solve the joint problem. This iterative
method proceeds in rounds in which each PM adjusts their trade lists.  In each
round each PM communicates their proposed trade list to a central
entity within the firm.  The central entity nets the trade lists and
broadcasts back to the PMs information about the net trades.  Based on this
information each PM modifies their optimization problem slightly, adding a
discount or premium to each asset, to account for the net trade transaction
cost.

The specific methods that the central entity and the PMs use come directly
from well-known methods for distributed optimization. Because of this, if the
adjustment rounds are continued, the PM trade lists converge to their joint
optimal values. But in this distributed method the PMs never divulge their
trade lists, or their strategies, to other PMs or indeed even the central
entity. In each round they simply take in information from the central
entity, modify their problems, and solve them to obtain the modified trade
list.

We observe that in practice, just a few rounds of adjustment are enough to
substantially reduce the transaction cost to the firm.

\subsection{Prior work}\label{s-lit}
\paragraph{Transaction costs.}
Transaction costs have been extensively studied in the context of portfolio
optimization. Many works have examined modeling transaction costs, and some
common models include a linear transaction cost model
\cite{li2019transaction,leland2000optimal}, a quadratic
transaction cost model \cite{chen2019note}, and a 3/2-power transaction cost
model \cite{boyd2017multi,boyd2024markowitz}. Other works have considered
piecewise-affine transaction cost models \cite{beraldi2019dealing} and a model
where the transaction cost is a fixed fraction of the portfolio value
\cite{morton1995optimal,mansini2015portfolio}. These works largely build on
Markowitz's mean-variance optimization framework \cite{markowitz1952portfolio}
and integrate transaction costs into a Markowitz portfolio optimization
problem. All of these papers, however, only consider the problem of optimizing
the portfolio of a single PM.

\paragraph{Distributed optimization.}
In the context of jointly optimizing the objectives of multiple agents in a
distributed setting, operator splitting methods such as the alternating
direction method of multipliers (ADMM) are well studied
\cite{boyd2010distributed,ryu2022large,chang2015multi,song2016fast,yang2022survey}.
We refer to a survey of methods in distributed optimization for a thorough
discussion of other methods used to solve multi-agent sharing problems
\cite{yang2019survey}. In portfolio optimization, distributed methods have been
studied to allow multiple agents to collaboratively optimize a single portfolio
based on observations each agent has made of asset returns and covariances
\cite{chen2013distributed}.

\section{Mathematical formulation}\label{s-math}
\subsection{Preliminaries}\label{s-prelim}
There are $M$ PMs at the same firm, each of whom manages a
portfolio of $N$ assets, but with different goals and mandates. Each portfolio
manager $i$ has a net asset value (NAV) $V^i > 0$, and the firm's total NAV
is $V = \sum_{i=1}^M V^i$. We define the relative
NAV weight of each PM as $\lambda^i = V^i / V$.

Each PM determines their trades in terms of portfolio weights.
We denote the trade weight vector of PM $i$ as $x^i \in
\reals^N$, where $x^i_j$ is the change in weight allocated to asset $j$, with
$x^i_j < 0$ meaning a reduction in weight. The corresponding trade list in
shares is $x^i V^i / p$, where $p \in \reals_{++}^N$ is the vector of asset
prices and the division is elementwise.

The firm collects the trades and executes the net trade in the market. Since
portfolios have different sizes, the firm's net trade weight is the
NAV-weighted sum of individual trade weights:
\[
z = \sum_{i=1}^M \lambda^i x^i.
\]
If $z_j = 0$, it means that shares of asset $j$ are exchanged between the $M$
funds in proportion to their NAVs, and none are transacted in the market. The
net trade in shares is $z V / p$.

When the trade lists are produced, the asset prices are given by $p_\text{ref}
\in \reals_{++}^N$, so the reference cost of the net trade is $(p_\text{ref})^T
(z V / p_\text{ref}) = V \ones^T z$. The trades
are executed at the realized prices $p_\text{real} \in \reals^N_{++}$, and the
realized cost differs due to market impact. The difference is interpreted as a
transaction cost, which is typically nonnegative, although it can be negative
if the trades of portfolio $i$ tend to go against the net trades, \ie, $x^i_j$
and $z_j$ have different signs.

\paragraph{Transaction cost models.}
There are many models of transaction cost. Since transaction costs depend on
the actual shares traded, we express them in terms of the net trade weight $z$
and the firm's total NAV. The simplest model includes only the bid-ask price
spread of each asset. Here, the reference price is the midpoint of the bid and
ask prices when the trade lists are produced, but we execute purchases at the
ask price and sales at the bid price. The resulting transaction cost is $(1/2)
\kappa^T |z|$, where $\kappa_j$ is the bid-ask spread of asset $j$ expressed as
a fraction of price (unitless), and the absolute value is taken elementwise.
More complex models also take into account the effect of large orders eating
through the order book. One common form, expressed in terms of trade weights,
is
\BEQ\label{e-3/2}
\phi^\text{tc}(z) = (1/2) \kappa_{\text{spread}}^T |z| + \kappa_{\text{impact}}^T |z|^{3/2},
\EEQ
with
\[
(\kappa_{\text{impact}})_j = \frac{b_j \nu_j}{\sqrt{\omega_j / V}}
\]
where $(\kappa_{\text{spread}})_j$ is the bid-ask spread of asset $j$ expressed
as a fraction of price (unitless); $b_j$ is the market impact coefficient for
asset $j$ (unitless); $\nu_j$ is the returns volatility of asset $j$ (unitless)
over the period considered for optimization; and $\omega_j$ is the dollar volume of
asset $j$ over the period considered for optimization. The ratio $\omega_j /
V$ can be interpreted as the market volume expressed in units of
portfolio weight. The $3/2$ power is taken elementwise. This models the
(predicted) transaction cost to execute the net trade $z$, expressed as a
fraction of NAV.

\paragraph{Portfolio manager objectives.} We suppose that PM $i$
has a closed, convex, proper objective function
$f^i:\reals^N\to\reals\cup\{+\infty\}$ such that $f^i(x)$ represents the
negative of the PM's expected return net of fees and costs
(possibly risk adjusted), expressed as a fraction of NAV, for trade weights $x\in
\reals^N$, for $i=1,\ldots,M$. This formulation ensures that all portfolio
manager objectives are of similar scale and directly comparable.
Additionally, we assume that portfolio
manager $i$ determines their trade weights $x^i$ such that
$x^i=\arg\min_{x}f^i(x)$. Note that $f^i$ can incorporate constraints. For
example, if PM $i$ solves the problem
\begin{equation}\label{e-generic-pm-problem}
\begin{array}{ll}
\mbox{minimize} & \tilde{f}^i(x^i) \\
\mbox{subject to} & x^i \in \mathcal{X}^i,
\end{array}
\end{equation}
with variable $x^i$, where $\mathcal{X}^i$ is the set of feasible trade weights
for PM $i$, then $f^i(x^i) = \tilde{f}^i(x^i) + I_{\mathcal{X}^i}(x^i)$, where
$I_{\mathcal{X}^i}$ is the indicator function of the feasible set
$\mathcal{X}^i$. This construction is used in \S\ref{s-backtest}. Common
constraints that may be included in $\mathcal{X}^i$ include risk limits
(bounding portfolio volatility), turnover limits (restricting the magnitude of
trades), leverage constraints (bounding the sum of absolute weights), and
active risk constraints (requiring that the optimized portfolio has limited
deviation from a target or benchmark portfolio).

\paragraph{Firm objective.} The firm's goal is to enable the portfolio
managers to achieve their objectives, while simultaneously minimizing the
transaction costs incurred by the firm as a whole. The firm's total \textit{ex
ante} penalty is therefore the NAV-weighted sum of the PMs' objective functions
plus the modeled transaction cost of the net trade. The problem of minimizing
this total penalty can be formulated as
\BEQ\label{e-prob-coop}
\begin{array}{ll}
\mbox{minimize}   & \sum_{i=1}^M \lambda^i f^i(x^i) + \gamma_\text{tc}\phi^\text{tc}(z) \\
\mbox{subject to} & z=\sum_{i=1}^M \lambda^i x^i
\end{array}
\EEQ
where $\gamma_\text{tc} > 0$ is a scaling parameter for the transaction cost
term. The NAV weights $\lambda^i$ ensure that the contribution of each PM to
the firm objective is proportional to their portfolio size. The inclusion of
$\gamma_\text{tc}$ is a practical matter: scaling the transaction cost term
often improves firm performance.

\paragraph{Note.} The formulation \eqref{e-prob-coop} describes the
problem of minimizing the penalty to the firm as a whole (or equivalently,
maximizing the firm's reward). Since the trade weights of the PMs
are coupled through the transaction cost term in the objective, the trades
of each PM will be influenced by the trades of the other
PMs. It is indeterminate how this will affect the performance
of every PM. We will show in appendix \ref{s-pm_results} that it
is demonstrably not the case that every PM will realize greater
performance by cooperating in this manner as opposed to acting independently. As
such, solving \eqref{e-prob-coop} is primarily of interest to firms that view
each PM as an investible asset. Such firms are only concerned
with the performance of the whole portfolio (the firm) and not the individual
assets (the PMs).

In principle, PMs could game the system by artificially scaling
up their objective functions (say, by a factor of 1000) to receive priority in
the firm optimization. However, we assume such gaming does not occur, as we
consider the cooperative case where PMs act in good faith toward
the firm objective.

\subsection{Alternating direction method of multipliers}\label{s-admm}
The alternating direction method of multipliers (ADMM) is a well-known
iterative method \cite{boyd2010distributed} that can be used to solve problems
of the form
\BEQ\label{e-prob-admm}
\begin{array}{ll}
\mbox{minimize}   & f(x) + g(z) \\
\mbox{subject to} & Ax+Bz=c
\end{array}
\EEQ
with variables $x\in\reals^n$ and $z\in\reals^m$, where $A\in\reals^{p\times
n}$, $B\in\reals^{p\times m}$, and $c\in\reals^p$ are constants. The ADMM
algorithm is given by
\BEAS
x^{k+1} &=& \argmin_x\left(f(x) + (u^k)^T (Ax+Bz^k-c) + \frac{\rho}{2}\|Ax+Bz^k-c\|_2^2\right) \\
z^{k+1} &=& \argmin_z\left(g(z) + (u^k)^T (Ax^{k+1}+Bz-c) + \frac{\rho}{2}\|Ax^{k+1}+Bz-c\|_2^2\right) \\
u^{k+1} &=& u^k + \varphi\rho(Ax^{k+1}+Bz^{k+1}-c)
\EEAS
where $\rho>0$ is a penalty parameter and $\varphi\in(0,\frac{1+\sqrt{5}}{2})$
is a step-size parameter \cite{ryu2022large}.

\paragraph{Reformulating the firm problem.}
To apply ADMM, we first introduce the NAV-scaled trade weights $\tilde{x}^i =
\lambda^i x^i$, which represent the contribution of PM $i$ to the firm's net
trade. The firm problem \eqref{e-prob-coop} can then be written as
\[
\begin{array}{ll}
\mbox{minimize} & \sum_{i=1}^M \lambda^i f^i(\tilde{x}^i/\lambda^i) + g(z) \\
\mbox{subject to} & \sum_{i=1}^M \tilde{x}^i - z = 0
\end{array}
\]
where $g(z) = \phi^\text{tc}(z) + I_{\mathcal{Z}}(z)$ is the sum of the
transaction cost and the indicator function of the feasible set $\mathcal{Z}$.
Following the dummy variables technique described in
\cite[Chapter~8]{ryu2022large}, we introduce new variables
$z^1,\ldots,z^M\in\reals^N$ and eliminate $z$ to obtain the equivalent problem
\[
\begin{array}{ll}
\mbox{minimize} & \sum_{i=1}^M \lambda^i f^i(\tilde{x}^i/\lambda^i) + g(\sum_{i=1}^M z^i) \\
\mbox{subject to} & \tilde{x}^i - z^i = 0, \quad i=1,\ldots,M.
\end{array}
\]
To improve practical convergence, we introduce a diagonal scaling matrix
$D\in\reals_{++}^{N\times N}$ and scale the constraints by $D$ to obtain the
equivalent problem
\BEQ\label{e-prob-admm-coop}
\begin{array}{ll}
\mbox{minimize} & \sum_{i=1}^M \lambda^i f^i(\tilde{x}^i/\lambda^i) + g(\sum_{i=1}^M z^i) \\
\mbox{subject to} & D\tilde{x}^i - Dz^i = 0, \quad i=1,\ldots,M.
\end{array}
\EEQ
As shown in appendix \ref{s-admm-appendix}, applying ADMM to problem
\eqref{e-prob-admm-coop} with initial dual variables satisfying
$u^{1,0}=\cdots=u^{M,0}$ allows us to derive the update rules. Writing the
updates in terms of the original PM trade weights $x^i = \tilde{x}^i/\lambda^i$
and defining the sharing update signal
\[
\ell^k = u^k+\frac{\rho}{M}\left(-Dz_\text{sum}^k+D\sum_{j=1}^M \lambda^j x^{j,k}\right),
\]
we obtain:
\BEAS
x^{i,k+1} &=& \argmin_x\left(\lambda^i f^i(x)
+ \lambda^i (\ell^k)^T Dx
+ \frac{\rho}{2}\|\lambda^i D(x-x^{i,k})\|_2^2\right) \\
z_\text{sum}^{k+1} &=& \argmin_z\left(g(z) - (u^k)^T Dz + \frac{\rho}{2M}\|Dz-{\textstyle\sum_{i=1}^M \lambda^i Dx^{i,k+1}}\|_2^2\right) \\
u^{k+1} &=& u^k + \frac{\varphi\rho}{M}(-Dz_\text{sum}^{k+1}+{\textstyle\sum_{j=1}^M \lambda^j Dx^{j,k+1}})
\EEAS

\subsection{Protocol}\label{s-algo}

The ADMM updates for solving \eqref{e-prob-admm-coop} can be computed in a
distributed manner, such that the PMs do not need to share their
objectives with each other or with the firm. Our distributed protocol works by
performing a fixed number of iterations of ADMM. While this does not produce an
exact solution to \eqref{e-prob-admm-coop}, we show empirically in
\S\ref{s-backtest} that even a few ADMM iterations is sufficient to
substantially refine the trade lists proposed by the PMs.

\begin{algorithm}[H]
\caption{Transaction cost mitigation protocol}
\label{a-algo}
\begin{algorithmic}[1]
\State \textbf{Initialization:}
\State \quad Fix iteration count $K$, step-size $\varphi$, penalty parameter $\rho$, and diagonal scaling matrix $D \in \reals_{++}^{N \times N}$
\State \quad Portfolio managers initialize: $x^{i,0} \gets \argmin_{x}f^i(x)$ for $i=1,\ldots,M$
\State \quad Central planner receives NAV-weighted net trade $\sum_{i=1}^M \lambda^i x^{i,0}$ and initializes $z_\text{sum}^0$, $u^0 \gets \mathbf{0}_N$
\State
\For{$k = 0, 1, \ldots, K-1$}
    \State \textbf{Step 1: Distributed update}
    \State \quad Central planner broadcasts signal:
    \State \quad\quad $\ell^{k} \gets u^k+\frac{\rho}{M}\left(-Dz_\text{sum}^k+D\sum_{i=1}^M \lambda^i x^{i,k}\right)$
    \State \quad Each PM $i$ updates:
    \State \quad\quad $x^{i,k+1} \gets \argmin_x\left(\lambda^i f^i(x) + \lambda^i (\ell^{k})^T Dx + \frac{\rho}{2}\|\lambda^i D(x-x^{i,k})\|_2^2\right)$
    \State
    \State \textbf{Step 2: Gathered update}
    \State \quad Central planner receives NAV-weighted net trade $\sum_{i=1}^M \lambda^i x^{i,k+1}$ and updates:
    \State \quad\quad $z_\text{sum}^{k+1} \gets \argmin_z\left(g(z) - (u^k)^T Dz + \frac{\rho}{2M}\left\|Dz-D\sum_{i=1}^M \lambda^i x^{i,k+1}\right\|_2^2\right)$
    \State \quad\quad $u^{k+1} \gets u^k + \frac{\varphi\rho}{M}\left(D\sum_{i=1}^M \lambda^i x^{i,k+1}-Dz_\text{sum}^{k+1}\right)$
\EndFor
\end{algorithmic}
\end{algorithm}

\paragraph{Comments.}
The PM update in Step 1 can be interpreted as re-solving the PM's original
optimization problem with two modifications: a linear price adjustment term
$(\ell^k)^T Dx$ that encodes discounts or premiums based on the net trade, and a
quadratic stability penalty $\frac{\rho}{2}\|\lambda^i D(x-x^{i,k})\|_2^2$ that
discourages large deviations from the previous trade list.

Notice that in the proposed algorithm, none of the PMs share their
trade lists with other PMs, and the PMs do not
share their objectives with the central planner. This structure preserves privacy
and prevents undesirable competition between PMs while also
minimizing costs for the entire firm.

\subsection{Hyperparameters and tuning}\label{s-hyper}

\paragraph{Choice of scaling matrix.}
The diagonal scaling matrix $D$ is used to improve
practical convergence of the algorithm.
The choice of $D$ can be motivated by approximating
the transaction cost function. We can approximate the $3/2$-power model with a
simpler one that is quadratic in trades, $\tilde{\phi}^\text{tc}(z) =
\kappa_\text{spread}^T|z| + z^T\diag(\kappa_\text{impact})z$. The Hessian of the
quadratic part of this approximation is $2\diag(\kappa_\text{impact})$. In ADMM,
the quadratic penalty on the consensus variable is weighted by $D^2$. A common
heuristic is to align this penalty with the curvature of the objective by setting
$D^2$ to be proportional to the Hessian of the centralized cost function's
quadratic part. This suggests the practical choice $D_{jj} =
\sqrt{2(\kappa_{\text{impact}})_j}$.

\paragraph{Choice of penalty and dual extrapolation parameters.}
From \cite{ryu2022large}, the objective value will converge to the global
optimum if $\rho>0$ and $\varphi\in(0,\frac{1+\sqrt{5}}{2})$. The penalty
parameter $\rho$ controls the strength of the consensus constraint, where
smaller values of $\rho$ result in larger changes in each iteration. While
larger steps may seem desirable, excessively large steps may lead to
overcorrection. The dual extrapolation parameter $\varphi$ controls the dual
variable update; a common choice is $\varphi=1$, though values up to
$\frac{1+\sqrt{5}}{2}\approx 1.618$ can accelerate convergence. Empirically, we
find the choices $\rho=10.0$ and $\varphi=1$ to work well in the backtest in
\S\ref{s-backtest}.

\paragraph{Choice of iteration count.}
As the iteration count $K$ is increased, the resulting trade lists will more
closely approximate the solution to \eqref{e-prob-coop}. As such, for firms
that prioritize the firm's returns and desire to solve \eqref{e-prob-coop}
exactly, a larger $K$ is preferred. For firms that prioritize the independence
of the PMs and wish only to provide a mechanism for the
PMs to account for transaction costs, a smaller value of $K$
(such as $K=1$) will minimally alter the initial trade lists of the portfolio
managers while potentially providing a substantial benefit in terms of
transaction cost mitigation.

\section{Extensions and variations}\label{s-extensions}

\paragraph{Firmwide constraints.} The firm might have constraints on the
aggregate trades and portfolio weights. For example, market availability
imposes box constraints on the net trade weights $z$. The firm might also have
constraints that limit the net exposure to any one asset class, limit the
leverage, or limit the turnover. We denote the set of feasible aggregate trade
weights as $\mathcal{Z}$. In this case, the firm's problem becomes
\[
\begin{array}{ll}
\mbox{minimize}   & \sum_{i=1}^M \lambda^i f^i(x^i) + \phi^\text{tc}(z) \\
\mbox{subject to} & z=\sum_{i=1}^M \lambda^i x^i\\
& z\in\mathcal{Z}.
\end{array}
\]
This problem can be written in the form \eqref{e-prob-coop} by taking
\[
\tilde{\phi}^\text{tc}(z) = \begin{cases}
\phi^\text{tc}(z), & \text{if $z\in\mathcal{Z}$} \\
+\infty, & \text{otherwise.}
\end{cases}
\]
This is a convex function when $\mathcal{Z}$ is convex. Algorithm~\ref{a-algo}
then applies.

\paragraph{Shorting costs.} When the firm nets positions across portfolio
managers, it only pays borrowing costs on the net short position rather than the
gross short positions of individual portfolios. This coupling can be
incorporated into the framework by adding a shorting cost term
$\phi^\text{short}(z)$ to the firm objective in \eqref{e-prob-coop}, and the
ADMM updates extend naturally. Other coupling terms arising from shared
resources or firmwide constraints can be handled similarly.

\section{Backtest}\label{s-backtest}
We present a backtest of our protocol
over historical data to illustrate the potential returns that can be
achieved by mitigating transaction costs in a firm-wide manner.
The code to reproduce our experiment can be found at
\begin{quote}
	\url{https://github.com/cvxgrp/coop_t_code}.
\end{quote}

\subsection{Backtest setup}\label{s-setup}

\paragraph{Data.}
We use historical market data for U.S.\ equities obtained from the LSEG (Refinitiv)
database. Our universe consists of $N = 434$ assets drawn from historical S\&P 500
constituents, filtered to those with sufficient data history and sorted by market
capitalization at the end of the sample period. For each asset, we obtain adjusted
daily trade prices, bid and ask prices, and trading volumes, where the adjustments
account for stock splits and dividends. The bid-ask spread, used in
the transaction cost model \eqref{e-3/2}, is computed as the absolute difference
between ask and bid prices. Missing values are forward-filled. We also obtain the
3-month U.S.\ Treasury Bill rate from the Federal Reserve Economic Data (FRED)
database, which serves as the risk-free rate $r_\text{rf}$ in the PM
objectives. The data spans July 2000 to April 2025.

This asset selection procedure introduces survivorship and lookahead bias, as we
select assets based on their market capitalization at the end of the sample period
rather than at each point in time. This is acceptable for our purposes: the goal is
to evaluate the relative benefit of the cooperative protocol compared to independent
optimization, not to demonstrate an investable strategy. Since both protocols
operate on the same asset universe, the comparison remains valid.

\paragraph{Transaction cost model.}
The transaction cost model \eqref{e-3/2} requires the spread parameter
$\kappa_\text{spread}$, the market impact coefficient $b$, the volatility
$\nu$, and the dollar volume $\omega$ for each asset. The spread parameter is
obtained directly from the bid-ask spread data. We set the market impact
coefficient $b_j = 1$ for all assets. The volatility $\nu_j$ is taken from the
diagonal of the factor model covariance estimate described below. The dollar
volume $\omega_j$ is computed as the product of the trading volume in shares and the
asset price. For shorting costs, we use the risk-free rate as a proxy for the
borrow cost, applied uniformly across all assets.

\paragraph{Policies.}
We simulate a firm with $M=4$ PMs. In each period of the backtest, each
PM determines their desired trades by solving the optimization problem
\begin{equation}\label{e-pm-problem}
\begin{array}{ll}
\mbox{minimize} & -\alpha^T w -r_\text{rf}c +\gamma_\text{risk}s_\text{risk}
+ \gamma_\text{turn}s_\text{turn} + \gamma_\text{tc}\phi^\text{tc}(z)
+ \gamma_\text{short}\phi^\text{short}(z) \\
\mbox{subject to} & w - w_\text{curr}=z\\
& 1-\mathbf{1}^T w = c\\
& \|w\|_1 \leq L\\
& |w| \leq C\\
&  \mathbf{1}^T \max(0, -w) \leq S\\
& \|w - w_\text{curr}\|_2 + |c - c_\text{curr}| \leq 2T + s_\text{turn}\\
& \|\Sigma^{1/2}w\|_2\leq \sigma_{\text{target}} + s_\text{risk}
\end{array}
\end{equation}
with $N$ tradable assets and:
\begin{itemize}
    \item Variables:
    \begin{itemize}
        \item $w \in \reals^N$, portfolio weights
        \item $c \in \reals$, cash position
        \item $z \in \reals^N$, trade weights
        \item $s_\text{risk} \in \reals_+$, risk slack variable
        \item $s_\text{turn} \in \reals_+$, turnover slack variable
    \end{itemize}
    \item Data:
    \begin{itemize}
        \item $w_\text{curr} \in \reals^N$, current portfolio weights
        \item $c_\text{curr} \in \reals$, current cash position
        \item $r_\text{rf} \in \reals$, risk-free rate
    \end{itemize}
    \item Models:
    \begin{itemize}
        \item $\alpha \in \reals^N$, expected asset returns
        \item $\Sigma \in \reals^{N \times N}$, asset covariance matrix
        \item $\phi^\text{tc}: \reals^N\to \reals_+$, transaction cost function
        \item $\phi^\text{short}: \reals^N\to \reals_{+}$, shorting cost function,
            $\phi^\text{short}(w) = r_\text{rf} \sum_{j=1}^N \max(0, -w_j)$
    \end{itemize}
    \item Parameters:
    \begin{itemize}
        \item $L \in \reals_+$, leverage limit
        \item $C \in \reals_+$, concentration limit
        \item $S \in \reals_+$, shorting limit
        \item $T \in \reals_+$, turnover limit
        \item $\sigma_{\text{target}} \in \reals_+$, target risk
        \item $\gamma_\text{risk} \in \reals_{++}$, risk penalty weight
        \item $\gamma_\text{turn} \in \reals_{++}$, turnover penalty weight
        \item $\gamma_\text{tc} \in \reals_{++}$, transaction cost penalty weight
        \item $\gamma_\text{short} \in \reals_{++}$, shorting cost penalty weight
    \end{itemize}
\end{itemize}

Each PM solves an optimization problem of the form above, with
identical constraint parameters but differing in their alpha estimates $\alpha$,
risk targets $\sigma_{\text{target}}$, tradable asset universes, and initial net
asset values (NAVs). The common parameters used across all PMs
are listed in Table~\ref{t-pm_params}. These values represent typical
institutional constraints and are consistent with parameters used in prior work
on systematic portfolio optimization \cite{boyd2024markowitz}. The risk targets
$\sigma_{\text{target}}$ are drawn uniformly at random from $[6\%, 15\%]$
annualized, and the initial NAVs are drawn log-uniformly from $10^{6.5}$ to
$10^{7.5}$ dollars. Each PM is assigned a random subset of 75\%
of the assets in the universe that they can trade and hold.

\begin{table}[H]
\centering
\begin{tabular}{lcc}
\hline
\textbf{Parameter} & \textbf{Symbol} & \textbf{Value} \\
\hline
Leverage limit & $L$ & 1.5 \\
Concentration limit & $C$ & 0.2 \\
Shorting limit & $S$ & 0.5 \\
Turnover limit & $T$ & 0.2 \\
Risk penalty weight & $\gamma_\text{risk}$ & 20 \\
Turnover penalty weight & $\gamma_\text{turn}$ & 1 \\
Transaction cost coefficient & $\gamma_\text{tc}$ & 0.15 \\
Shorting cost coefficient & $\gamma_\text{short}$ & 1 \\
\hline
\end{tabular}
\caption{Common constraint and penalty parameters used across all PMs.}
\label{t-pm_params}
\end{table}

\noindent All PMs
use a common risk model. The details of the PM alphas, asset
universes, and common risk model are described in the following paragraphs.

This formulation is inspired by the paper \cite{boyd2024markowitz}, and is intended
to represent a reasonable method for systematic portfolio optimization. This problem
is an instance of problem \eqref{e-generic-pm-problem}.

\paragraph{Cooperative protocols.} We implement and compare four models
(protocols) of PM cooperation:
\begin{itemize}
    \item \texttt{independent}: The baseline protocol where each portfolio
    manager solves their optimization problem independently, accounting for
    their own estimated transaction costs without coordination.
    \item \texttt{full\_cooperative}: An idealized protocol that solves problem
    \eqref{e-prob-coop} directly, combining all PMs' objectives
    into a single optimization problem with a joint transaction cost term applied
    to the net firm trade and a firm-wide shorting cost term as described in
    \S\ref{s-extensions}. The transaction and shorting cost terms are scaled by
    $\gamma_\text{tc}$ and $\gamma_\text{short}$, respectively, matching the
    parameters used in the individual PM policies. This requires gathering all
    the problem data from the PMs to solve as implemented.
    \item \texttt{admm\_2\_iter}: The distributed ADMM protocol described in
    \S\ref{s-algo} with $K=2$ iterations, which approximately solves the same
    problem as \texttt{full\_cooperative} without requiring PMs
    to share their complete objectives with the firm or each other. The ADMM
    hyperparameters are set as described in \S\ref{s-hyper}.
    \item \texttt{admm\_5\_iter}: The ADMM protocol with $K=5$ iterations,
    providing a closer approximation to the \texttt{full\_cooperative} solution.
\end{itemize}

\paragraph{Synthetic alpha generation.}
To evaluate our protocol under controlled conditions, we generate synthetic
alpha estimates with known statistical properties rather than using proprietary
trading signals. Each PM receives a noised estimate of future
returns, calibrated to achieve a target Information Coefficient (IC)---the
cross-sectional correlation between the alpha estimate and realized returns.

For PM $i$, the synthetic alpha estimate at time $t$ is
\[
\hat{R}_t^{(i)} = c^{(i)}(R_t + E_t^{(i)}),
\]
where $R_t \in \reals^N$ is the vector of realized returns over a 42-day
forward-looking window, $E_t^{(i)} \in \reals^N$ is a noise vector, and
$c^{(i)} > 0$ is a scaling factor. The noise variance is calibrated such that
the correlation between $\hat{R}_t^{(i)}$ and $R_t$ equals the target IC for
strategy $i$. Specifically, if the target IC is $\rho^{(i)}$, then the scaling
factor is $c^{(i)} = (\rho^{(i)})^2$ and the noise scaling factor is
$v^{(i)} = \sqrt{1/(\rho^{(i)})^2 - 1}$.

The noise process $E_t$ (stacked across all strategies) follows a vector
autoregressive process of order 1 (VAR(1)):
\[
E_t = \Phi E_{t-1} + U_t,
\]
where $\Phi$ is a block-diagonal matrix with strategy-specific temporal
autocorrelation coefficients $\phi^{(i)}$ on the diagonal blocks, and $U_t$ is
an innovation vector. The stationary covariance of $E_t$ has the Kronecker
structure $\Sigma_E = S_E \otimes \Sigma_{\text{Asset}}$, where
$\Sigma_{\text{Asset}}$ is the empirical asset return covariance and $S_E$
encodes cross-strategy correlations scaled by the noise factors $v^{(i)}$. The
innovation covariance $\Sigma_U$ is determined from the discrete-time Lyapunov
equation $\Sigma_E = \Phi \Sigma_E \Phi^T + \Sigma_U$. The synthetic alpha
parameters are summarized in Table~\ref{t-alpha_params}.

\begin{table}[H]
\centering
\begin{tabular}{lc}
\hline
\textbf{Parameter} & \textbf{Value} \\
\hline
Target IC (per PM) & $\mathcal{U}[0.06, 0.10]$ \\
Temporal autocorrelation $\phi^{(k)}$ & $\mathcal{U}[0.75, 0.85]$ \\
Cross-strategy noise correlation & $0.3$ \\
Forward return horizon & 42 days \\
\hline
\end{tabular}
\caption{Synthetic alpha generation parameters. $\mathcal{U}[a,b]$ denotes a
uniform distribution over $[a,b]$.}
\label{t-alpha_params}
\end{table}

Note that these alpha estimates are not
investable in practice, as they rely on noised versions of true future returns,
which are unknowable. The PMs and their alpha models should be
treated as proxies for realistic systematic portfolio optimization strategies
rather than as implementable trading strategies.

\paragraph{Covariance matrix estimation.} All PMs share a common
risk model based on a low-rank factor structure. For each time $t$, we estimate
the covariance matrix $\Sigma_t = F_t F_t^T + D_t$, where $F_t \in \reals^{N
\times J}$ is the factor loading matrix with $J = 15$ factors, and $D_t \in
\reals^{N \times N}$ is a diagonal matrix of idiosyncratic variances.

The estimation procedure uses a 42-day forward-looking window. At each time $t$,
we compute the realized returns over the subsequent 42 trading days and estimate
per-asset volatilities as the standard deviation over this window. Returns are
standardized by these volatilities and winsorized at $\pm 4.2$ standard
deviations to reduce the influence of outliers. We then compute the sample
correlation matrix of the winsorized standardized returns and extract the top
$J$ eigenvectors and eigenvalues. The factor loadings are constructed as $F_t =
\text{diag}(\sigma_t) Q_t \Lambda_t^{1/2}$, where $\sigma_t$ is the vector of
asset volatilities, $Q_t$ contains the top $J$ eigenvectors, and $\Lambda_t$ is
the diagonal matrix of corresponding eigenvalues. The idiosyncratic variances
are set to the residual variance not explained by the factors.

This forward-looking estimation is not implementable in practice, as it uses
future information. However, for the purpose of evaluating the cooperative
trading protocol under controlled conditions, this approach ensures that the
risk model accurately reflects the true covariance structure during the backtest
period.

\subsection{Results}\label{s-results}
In Table~\ref{t-performance_stats}, we present the top-level performance statistics
for the firm in each of the backtest scenarios. Since the underlying strategies are
synthetic, the absolute performance figures are not meaningful in themselves; what
matters is the relative improvement achieved through cooperation. We observe that
the cooperative protocols all significantly outperform the independent protocol in
terms of return, volatility, and Sharpe ratio. In comparison, the difference between the ADMM protocols
and the full joint protocol is not as pronounced, and all three attain nearly identical
Sharpe ratios. We break down the performance by PM in appendix \ref{s-pm_results}.

\begin{table}[H]
\centering
\begin{tabular}{lccc}
\hline
    & \textbf{Return} & \textbf{Volatility} & \textbf{Sharpe} \\
\hline
\texttt{independent} & $15.63\%$ & $9.77\%$ & $1.60$ \\
\texttt{full\_cooperative} & $17.59\%$ & $8.58\%$ & $2.05$ \\
\texttt{admm\_2\_iter} & $18.70\%$ & $9.16\%$ & $2.04$ \\
\texttt{admm\_5\_iter} & $18.63\%$ & $8.96\%$ & $2.08$ \\
\hline
\end{tabular}
\caption{Performance statistics by backtest scenario.}
\label{t-performance_stats}
\end{table}

\begin{figure}[H]
\centering
\includegraphics[width=0.9\textwidth]{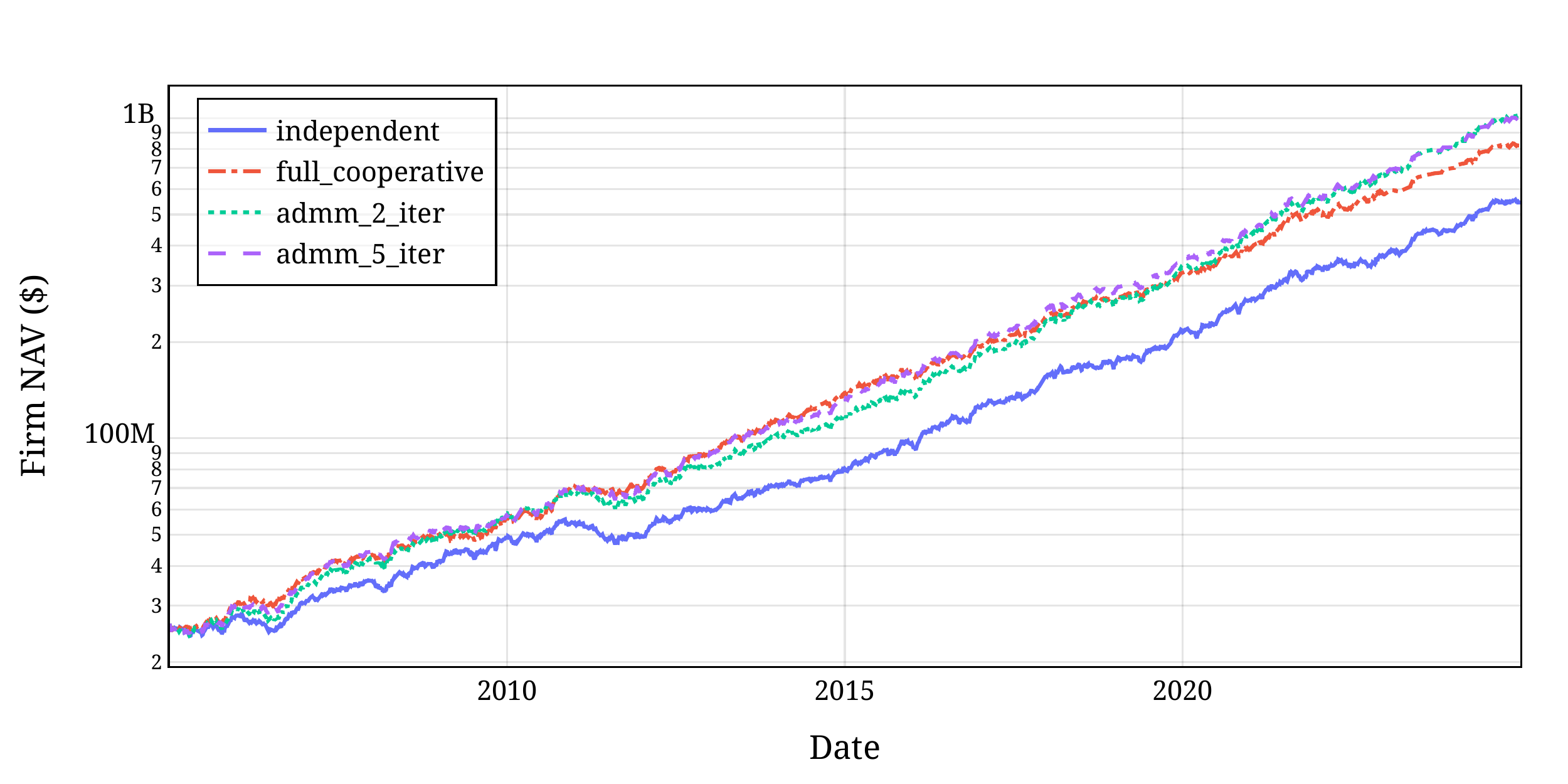}
\caption{Firm's cumulative returns by backtest scenario.}
\label{f-returns}
\end{figure}

\begin{figure}[H]
\centering
\includegraphics[width=0.9\textwidth]{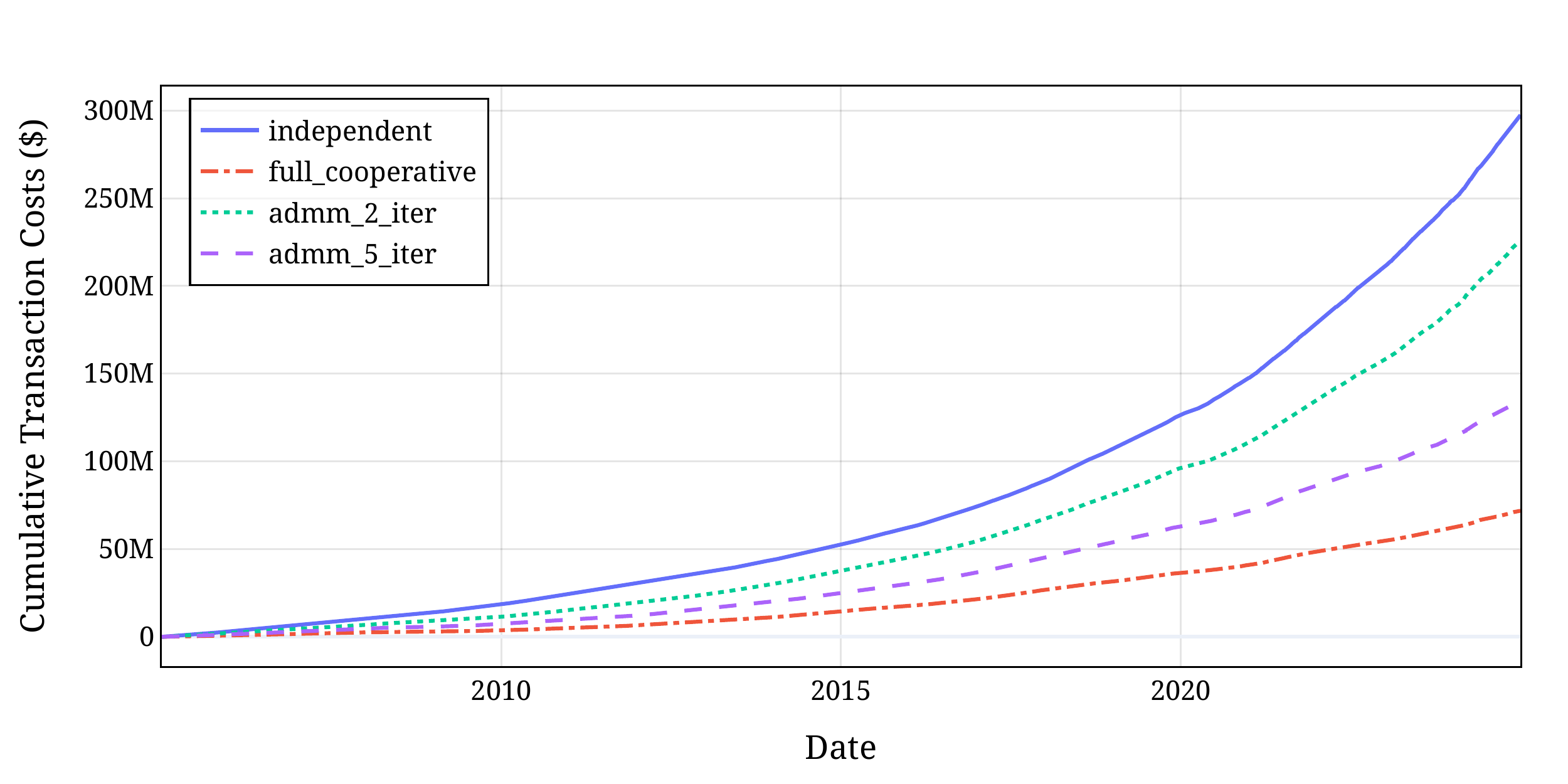}
\caption{Firm's cumulative transaction costs by backtest scenario.}
\label{f-tcosts}
\end{figure}

We show the firm's cumulative returns and cumulative transaction costs in each backtest
scenario in Figures~\ref{f-returns} and~\ref{f-tcosts}, respectively.
As should be expected, the full joint protocol achieves the lowest
transaction costs, and more iterations of the ADMM protocol seem to result in lower
transaction costs. Performing 2 iterations of ADMM already achieves a substantial
reduction in transaction costs, and performing 5 iterations yields approximately 75\%
of the savings of the full joint protocol, and there is very little difference between
the cumulative returns of the ADMM protocols and the full joint protocol.

\section{Conclusion}\label{s-conclusion}
We have presented a distributed protocol, based on ADMM, that enables portfolio
managers within a firm to coordinate their trades and reduce transaction costs
without revealing their objectives to each other or to the firm. Our
backtest demonstrates that even a small number of coordination rounds can
substantially reduce the transaction costs incurred by the firm on its net
trades.

It is important to note the limitations of this approach. The protocol does not
guarantee improved returns or Sharpe ratios for the firm or for individual
PMs. The effect on PM-level performance is
indeterminate: while we expect every PM to experience reduced
transaction costs through cooperation, this does not necessarily translate to
improved returns, as the coordination may alter the composition of trades in
ways that affect other aspects of performance. What this method does provide is
a principled mechanism for reducing frictional costs that would otherwise erode
returns.

For firms structured as collections of independent sleeves or subsidiary funds,
where net trades are executed centrally, this protocol offers a practical path
to capturing the benefits of trade netting while preserving the autonomy and
privacy of individual PMs.

\section{Acknowledgements}
We thank Ronald Kahn for helpful discussions about transaction cost models.

\clearpage
\bibliographystyle{alpha}
\bibliography{refs}

\clearpage
\appendix
\section{ADMM for the firm problem}\label{s-admm-appendix}
Recall that we define the NAV-scaled trade weights $\tilde{x}^i = \lambda^i
x^i$. Applying ADMM to \eqref{e-prob-admm-coop} yields the iterates (written in
terms of $\tilde{x}^i$):
\BEAS
\tilde{x}^{i,k+1} &=& \argmin_{\tilde{x}}\left(\lambda^i f^i(\tilde{x}/\lambda^i) + (u^{i,k})^T (D\tilde{x}-Dz^{i,k}) + \frac{\rho}{2}\|D\tilde{x}-Dz^{i,k}\|_2^2\right) \\
z^{k+1} &=& \argmin_{(z^1,\ldots,z^M)}
\left(g({\textstyle\sum_{i=1}^M z^i}) + \sum_{i=1}^M\left[(u^{i,k})^T
(D\tilde{x}^{i,k+1}-Dz^i) + \frac{\rho}{2}\|D\tilde{x}^{i,k+1}-Dz^i\|_2^2\right]
\right) \\
u^{i,k+1} &=& u^{i,k} + \varphi\rho(D\tilde{x}^{i,k+1}-Dz^{i,k+1}).
\EEAS
If ADMM is initialized with $u^{1,0}=\cdots=u^{M,0}$, then
$u^{1,k}=\cdots=u^{M,k}$ and $D\tilde{x}^{1,k}-Dz^{1,k}=\cdots=D\tilde{x}^{M,k}-Dz^{M,k}$ at
each iteration $k$. We can show this inductively. Suppose
$u^{1,k}=\cdots=u^{M,k}$ and rewrite the $z$-update by introducing auxiliary
variable $z_\text{sum}$, solving the problem
\[
\begin{array}{ll}
\mbox{minimize} & g(z_\text{sum}) - (u^k)^T(Dz_\text{sum}) + \sum_{i=1}^M \frac{\rho}{2}\|D\tilde{x}^{i,k+1}-Dz^i\|_2^2\\
\mbox{subject to} & \sum_{i=1}^M z^i = z_\text{sum}
\end{array}
\]
with variables $z^1,\ldots,z^M$ and $z_\text{sum}$. Solving first for
$z^1,\ldots,z^M$ and then for $z_\text{sum}$ gives
$z^{i,k+1}=\tilde{x}^{i,k+1}-\frac{1}{M}\sum_{j=1}^M
\tilde{x}^{j,k+1}+\frac{1}{M}z_\text{sum}^{k+1}$, where $z_\text{sum}^{k+1}$ is given
by
\[
z_\text{sum}^{k+1} = \argmin_z\left(g(z) - (u^k)^T Dz + \frac{\rho}{2M}\|Dz-{\textstyle\sum_{i=1}^M D\tilde{x}^{i,k+1}}\|_2^2\right).
\]
In particular,
$D\tilde{x}^{i,k+1}-Dz^{i,k+1}=-\frac{1}{M}Dz_\text{sum}^{k+1}+\frac{1}{M}\sum_{j=1}^M
D\tilde{x}^{j,k+1}$ for $i=1,\ldots,M$, which is constant across $i$. Consequently,
from the $u$-update
\[
u^{i,k+1} = u^k + \varphi\rho(D\tilde{x}^{i,k+1}-Dz^{i,k+1})
= u^k + \frac{\varphi\rho}{M}(-Dz_\text{sum}^{k+1}+{\textstyle\sum_{j=1}^M D\tilde{x}^{j,k+1}}),
\]
we also have $u^{1,k+1}=\cdots=u^{M,k+1}$.

We can also rewrite the $\tilde{x}$-update in terms of $z_\text{sum}$. Substituting
$z^{i,k}=\tilde{x}^{i,k}-\frac{1}{M}\sum_{j=1}^M \tilde{x}^{j,k}+\frac{1}{M}z_\text{sum}^k$
into the $\tilde{x}$-update, expanding the quadratic term, and eliminating terms
constant with respect to $\tilde{x}$, and defining the sharing update signal
\[
\ell^k = u^k+\frac{\rho}{M}\left(-Dz_\text{sum}^k+D\sum_{j=1}^M \tilde{x}^{j,k}\right),
\]
we get
\[
\tilde{x}^{i,k+1} = \argmin_{\tilde{x}}\left(\lambda^i f^i(\tilde{x}/\lambda^i)
+ (\ell^k)^T D\tilde{x}
+ \frac{\rho}{2}\|D\tilde{x}-D\tilde{x}^{i,k}\|_2^2\right).
\]
Substituting $\tilde{x} = \lambda^i x$ and $\tilde{x}^{j,k} = \lambda^j x^{j,k}$
and simplifying, we obtain the update in terms of the original PM trade weights $x^i$:
\[
x^{i,k+1} = \argmin_x\left(\lambda^i f^i(x)
+ \lambda^i (\ell^k)^T Dx
+ \frac{\rho}{2}\|\lambda^i D(x-x^{i,k})\|_2^2\right).
\]
Combining these new iterates gives the update rules described in
\S\ref{s-admm}.

\section{Portfolio manager results}\label{s-pm_results}
Figures~\ref{f-pm_returns} and~\ref{f-pm_tcosts} show the PM's
cumulative returns and cumulative transaction costs by backtest scenario, respectively.
Backtest statistics for each PM are presented in Table~\ref{t-performance_stats_pm}.
We do expect that every PM experiences
lower transaction costs, which we see in Figure~\ref{f-pm_tcosts}, but this method
offers no guarantee that every PM realizes greater returns or Sharpe
ratio through cooperation. In particular, we observe that PM 2 experiences worse
performance through cooperation compared to optimizing independently post 2010.

\begin{figure}[H]
\centering
\begin{subfigure}[b]{0.9\textwidth}
    \centering
    \includegraphics[width=\textwidth]{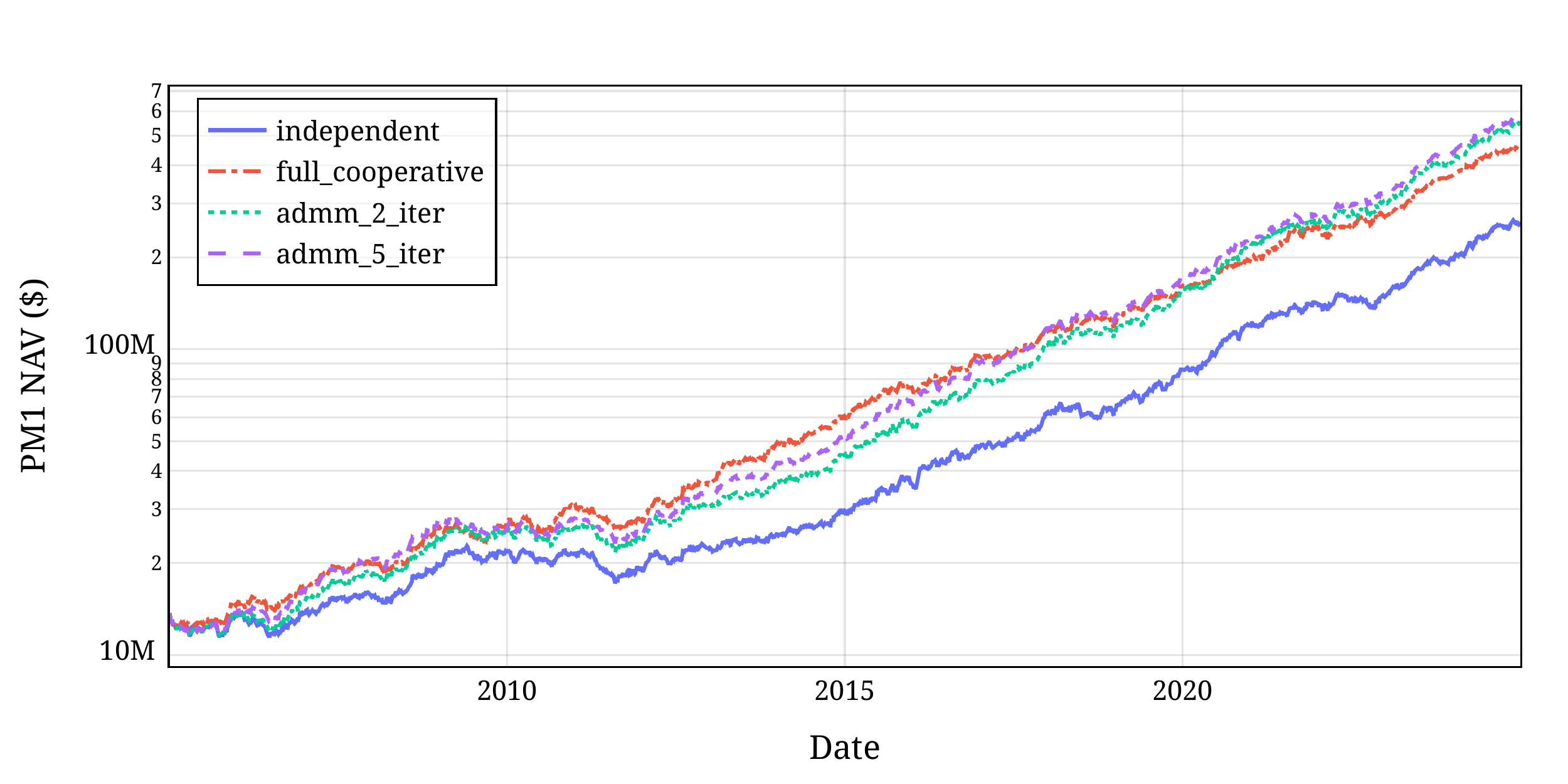}
    \caption{PM 1}
\end{subfigure}
\end{figure}

\begin{figure}[H]
\ContinuedFloat
\centering
\begin{subfigure}[b]{0.9\textwidth}
    \centering
    \includegraphics[width=\textwidth]{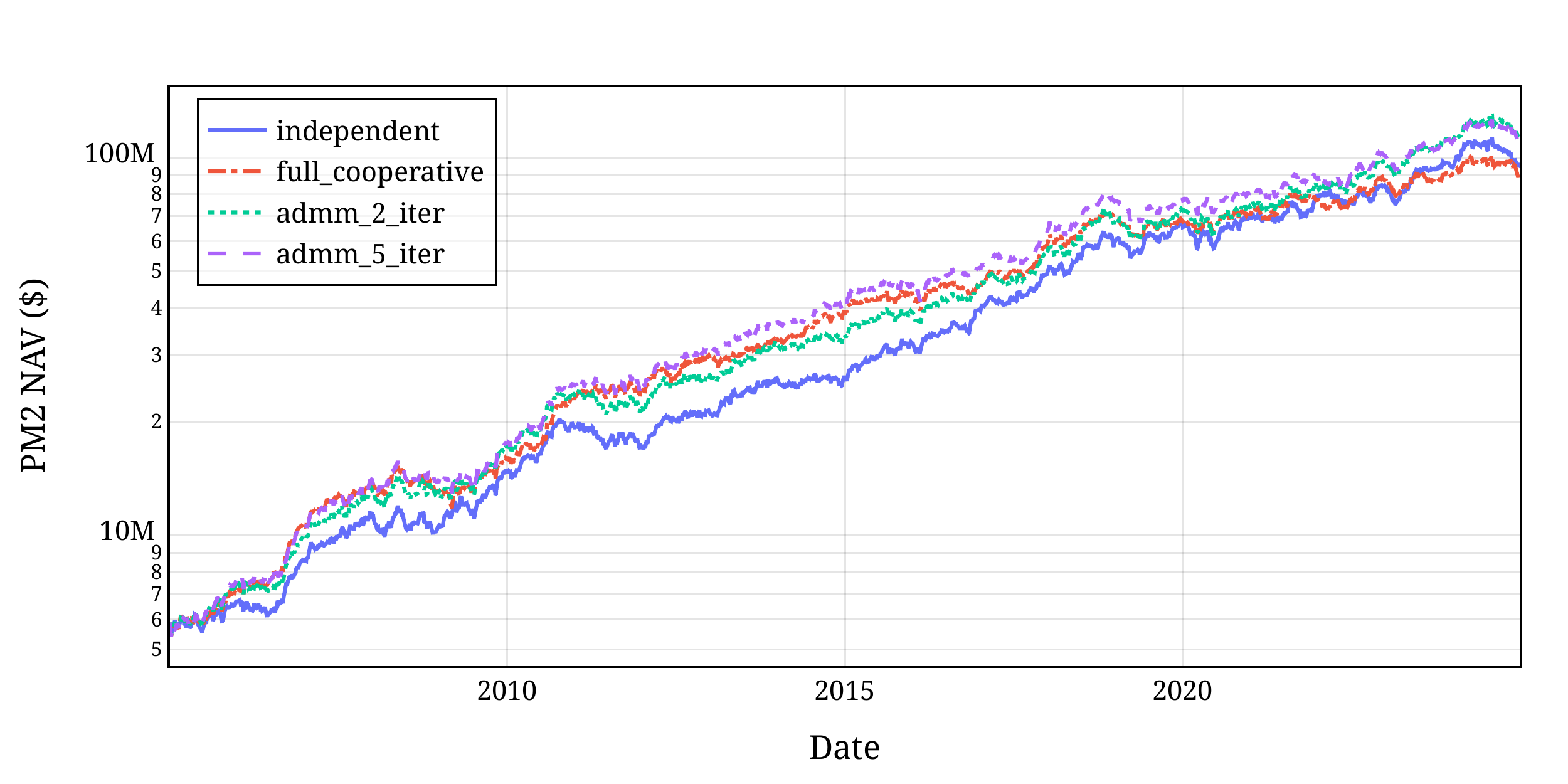}
    \caption{PM 2}
\end{subfigure}
\end{figure}

\begin{figure}[H]
\ContinuedFloat
\centering
\begin{subfigure}[b]{0.9\textwidth}
    \centering
    \includegraphics[width=\textwidth]{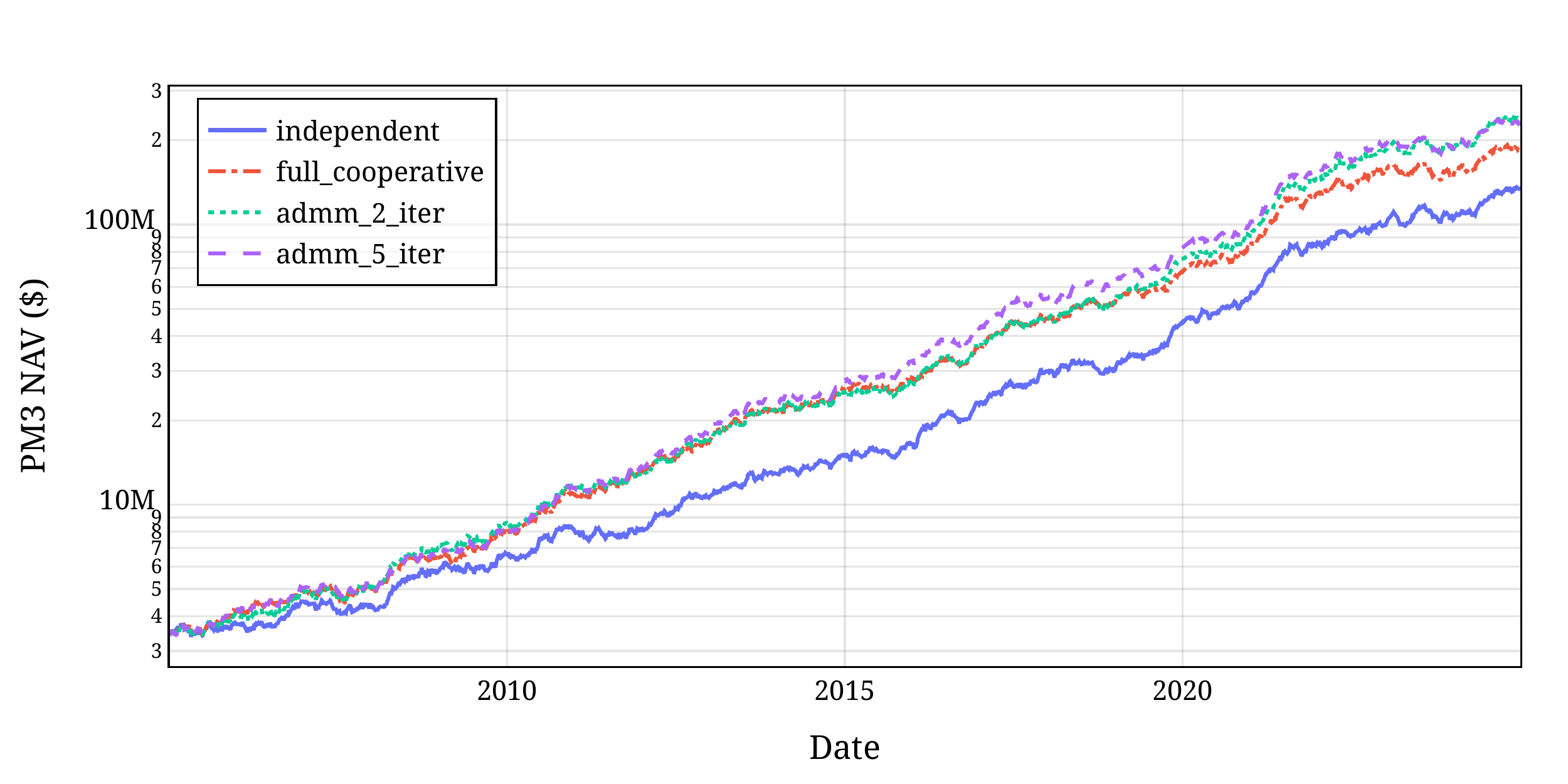}
    \caption{PM 3}
\end{subfigure}
\end{figure}

\begin{figure}[H]
\ContinuedFloat
\centering
\begin{subfigure}[b]{0.9\textwidth}
    \centering
    \includegraphics[width=\textwidth]{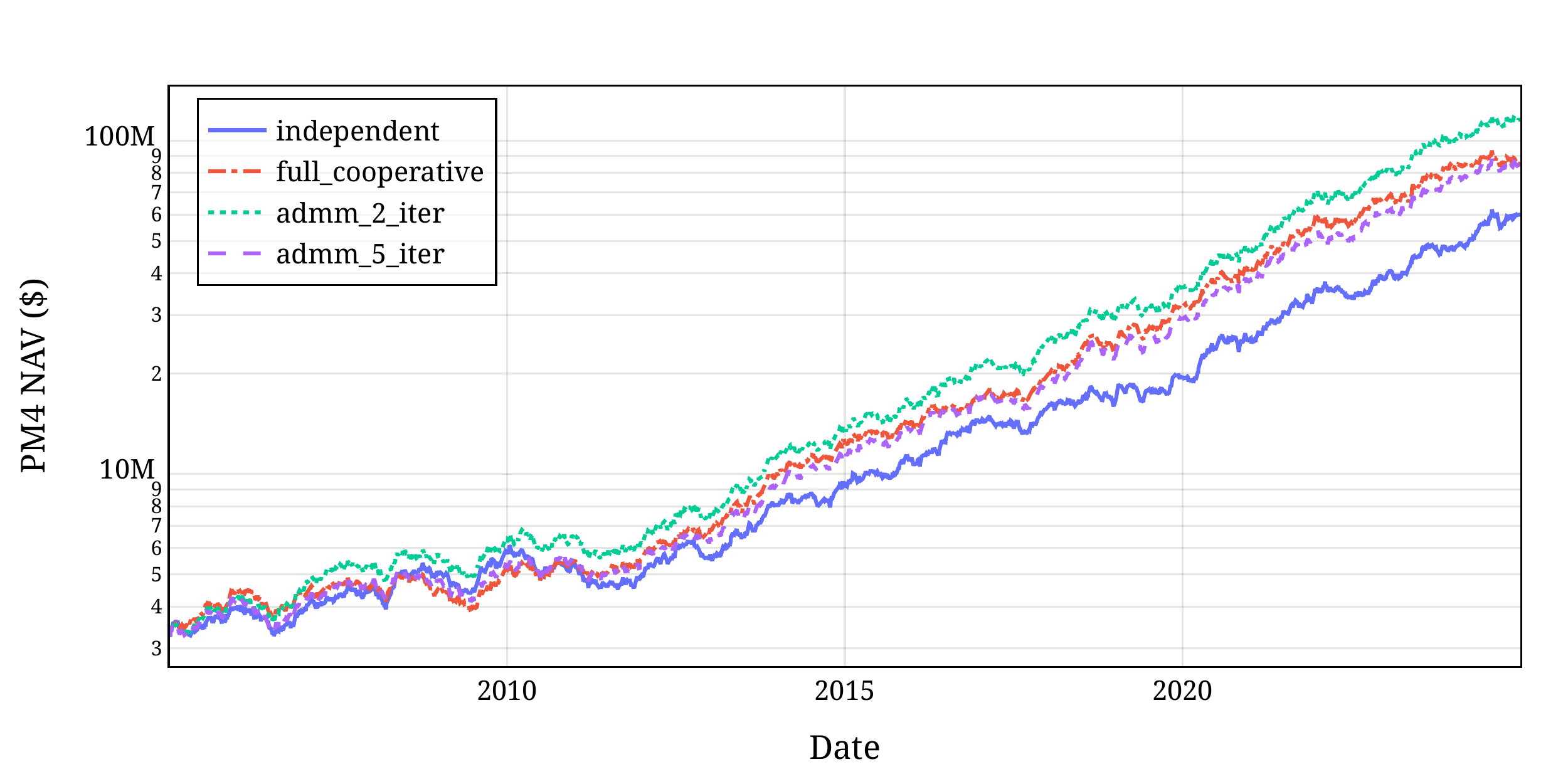}
    \caption{PM 4}
\end{subfigure}
\caption{Portfolio managers' cumulative returns by backtest scenario.}
\label{f-pm_returns}
\end{figure}

\begin{figure}[H]
\centering
\begin{subfigure}[b]{0.9\textwidth}
    \centering
    \includegraphics[width=\textwidth]{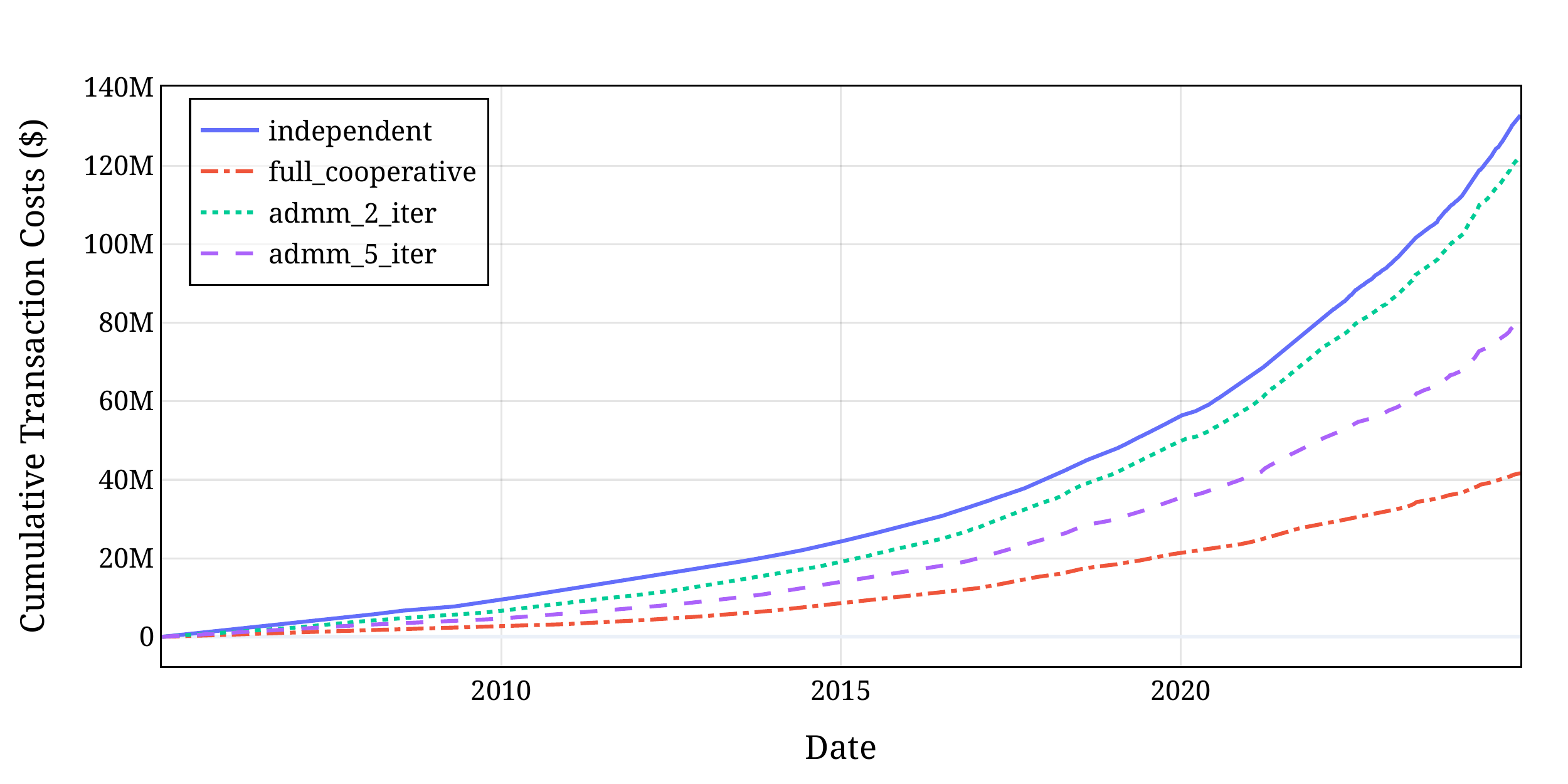}
    \caption{PM 1}
\end{subfigure}
\end{figure}

\begin{figure}[H]
\ContinuedFloat
\centering
\begin{subfigure}[b]{0.9\textwidth}
    \centering
    \includegraphics[width=\textwidth]{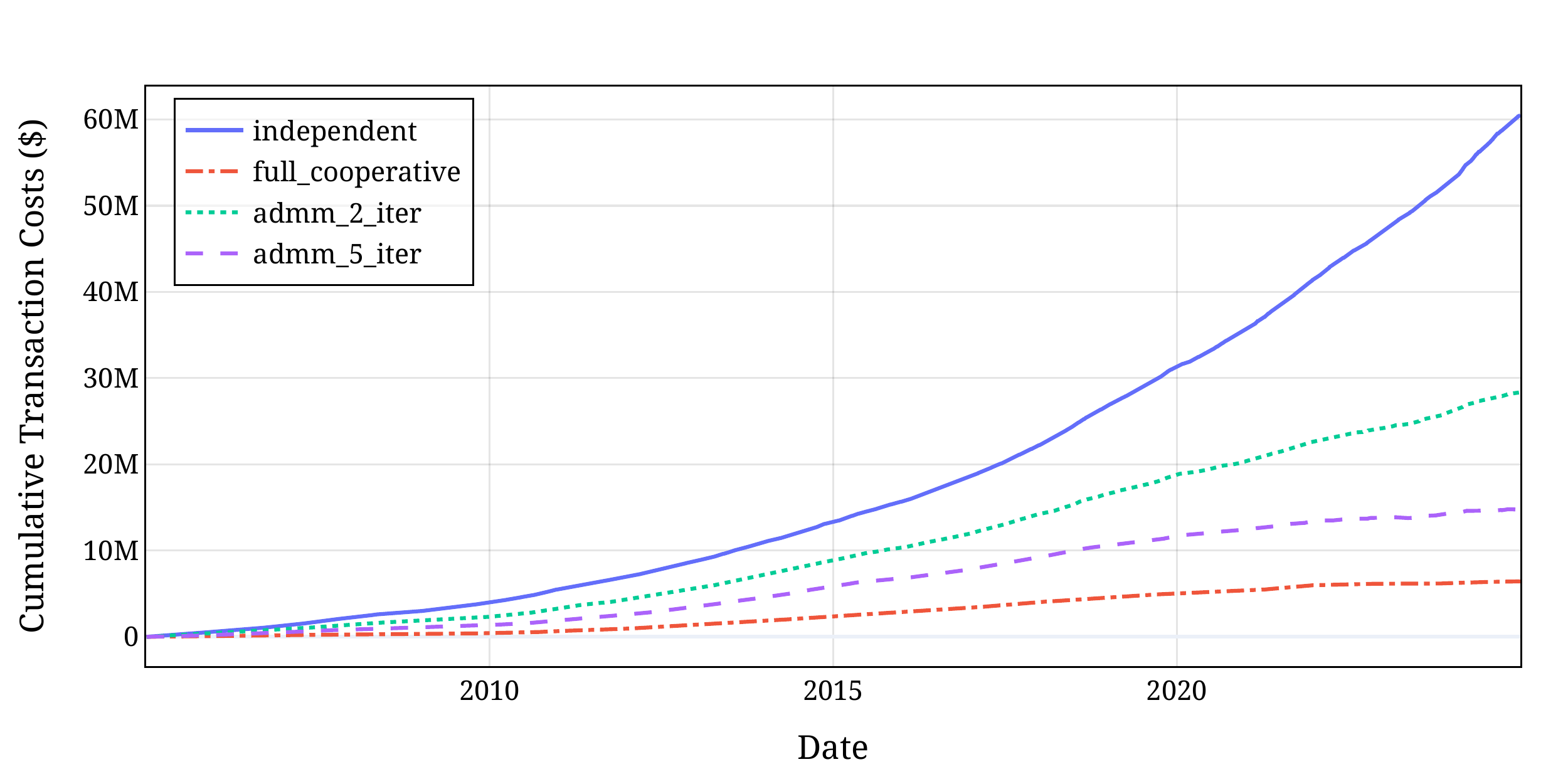}
    \caption{PM 2}
\end{subfigure}
\end{figure}

\begin{figure}[H]
\ContinuedFloat
\centering
\begin{subfigure}[b]{0.9\textwidth}
    \centering
    \includegraphics[width=\textwidth]{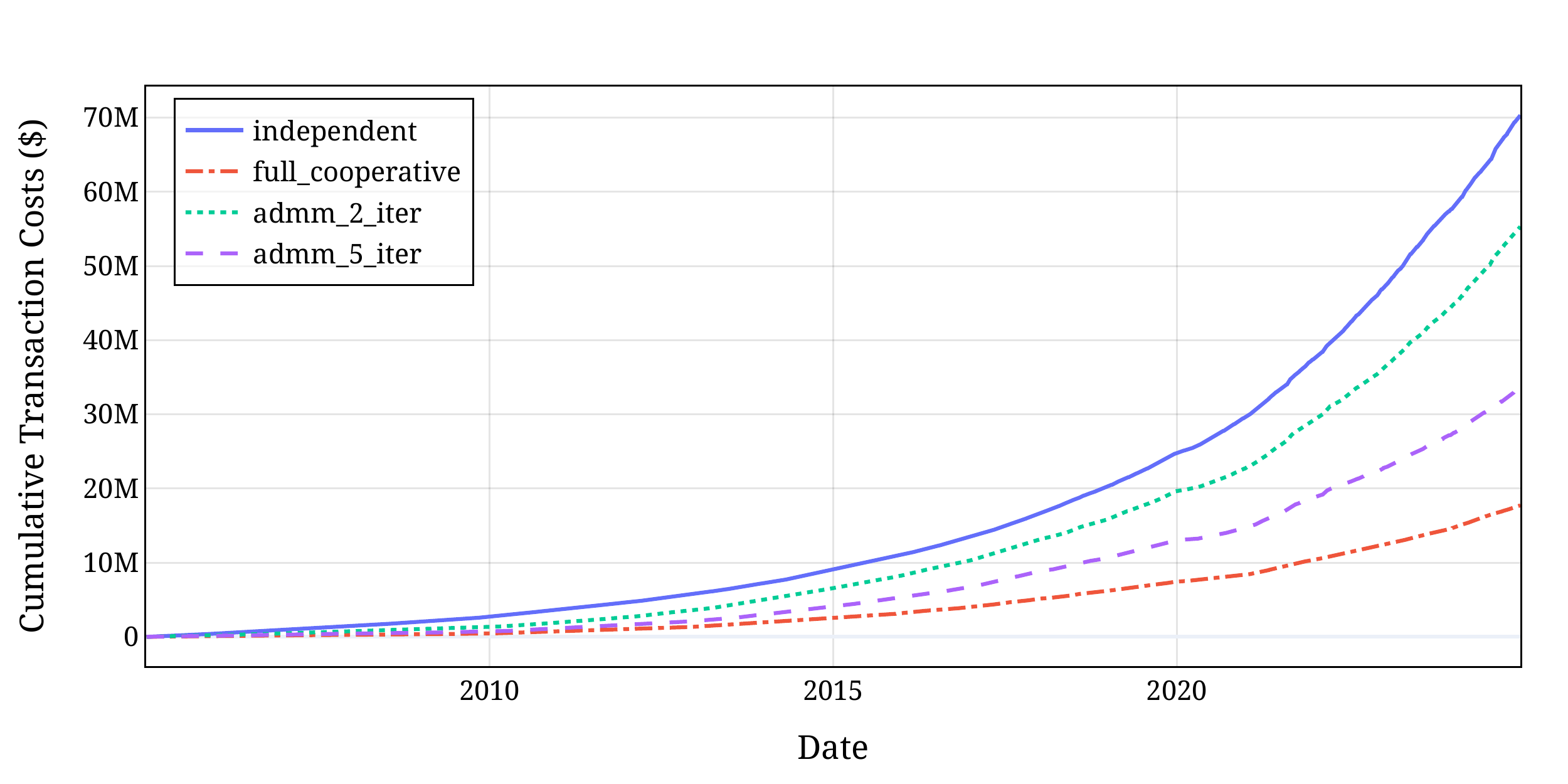}
    \caption{PM 3}
\end{subfigure}
\end{figure}

\begin{figure}[H]
\ContinuedFloat
\centering
\begin{subfigure}[b]{0.9\textwidth}
    \centering
    \includegraphics[width=\textwidth]{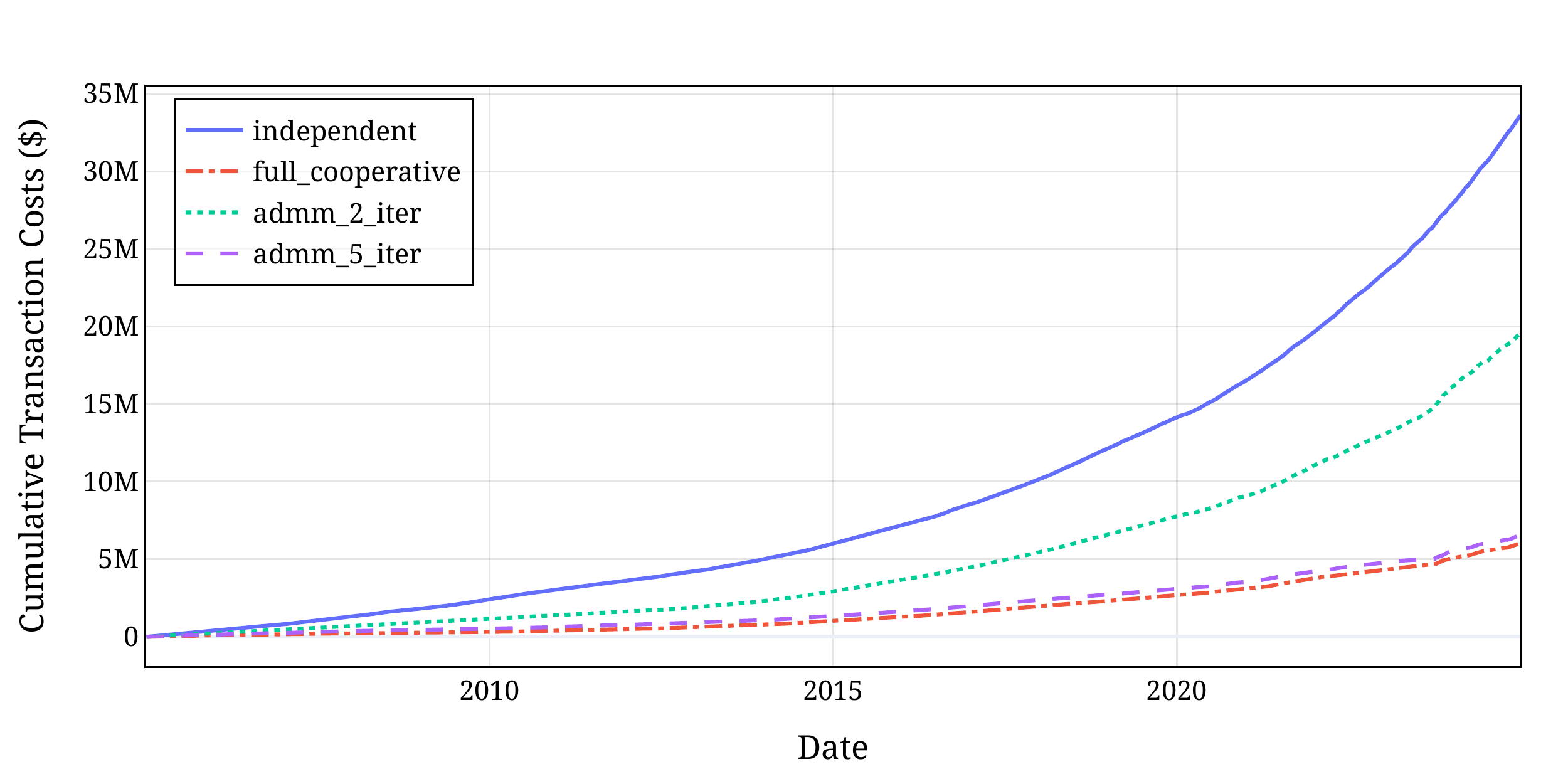}
    \caption{PM 4}
\end{subfigure}
\caption{Portfolio managers' cumulative transaction costs by backtest scenario.}
\label{f-pm_tcosts}
\end{figure}

\begin{table}[H]
\centering
\begin{subtable}[b]{0.9\textwidth}
\centering
\begin{tabular}{lccc}
\hline
& \textbf{Return} & \textbf{Volatility} & \textbf{Sharpe} \\
\hline
\texttt{independent} & $15.45\%$ & $12.57\%$ & $1.23$ \\
\texttt{full\_cooperative} & $18.23\%$ & $11.14\%$ & $1.64$ \\
\texttt{admm\_2\_iter} & $19.15\%$ & $11.85\%$ & $1.62$ \\
\texttt{admm\_5\_iter} & $19.42\%$ & $11.65\%$ & $1.67$ \\
\hline
\end{tabular}
\caption{PM 1}
\label{t-performance_stats_pm1}
\end{subtable}
\end{table}

\begin{table}[H]
\centering
\ContinuedFloat
\begin{subtable}[b]{0.9\textwidth}
\centering
\begin{tabular}{lccc}
\hline
& \textbf{Return} & \textbf{Volatility} & \textbf{Sharpe} \\
\hline
\texttt{independent} & $14.79\%$ & $13.38\%$ & $1.10$ \\
\texttt{full\_cooperative} & $14.39\%$ & $12.39\%$ & $1.16$ \\
\texttt{admm\_2\_iter} & $15.70\%$ & $12.95\%$ & $1.21$ \\
\texttt{admm\_5\_iter} & $15.49\%$ & $12.67\%$ & $1.22$ \\
\hline
\end{tabular}
\caption{PM 2}
\label{t-performance_stats_pm2}
\end{subtable}
\end{table}

\begin{table}[H]
\centering
\ContinuedFloat
\begin{subtable}[b]{0.9\textwidth}
\centering
\begin{tabular}{lccc}
\hline
& \textbf{Return} & \textbf{Volatility} & \textbf{Sharpe} \\
\hline
\texttt{independent} & $19.02\%$ & $12.34\%$ & $1.54$ \\
\texttt{full\_cooperative} & $20.54\%$ & $11.34\%$ & $1.81$ \\
\texttt{admm\_2\_iter} & $21.87\%$ & $11.69\%$ & $1.87$ \\
\texttt{admm\_5\_iter} & $21.65\%$ & $11.42\%$ & $1.90$ \\
\hline
\end{tabular}
\caption{PM 3}
\label{t-performance_stats_pm3}
\end{subtable}
\end{table}

\begin{table}[H]
\centering
\ContinuedFloat
\begin{subtable}[b]{0.9\textwidth}
\centering
\begin{tabular}{lccc}
\hline
& \textbf{Return} & \textbf{Volatility} & \textbf{Sharpe} \\
\hline
\texttt{independent} & $15.42\%$ & $13.09\%$ & $1.18$ \\
\texttt{full\_cooperative} & $17.10\%$ & $12.28\%$ & $1.39$ \\
\texttt{admm\_2\_iter} & $18.60\%$ & $12.28\%$ & $1.51$ \\
\texttt{admm\_5\_iter} & $17.06\%$ & $12.31\%$ & $1.39$ \\
\hline
\end{tabular}
\caption{PM 4}
\label{t-performance_stats_pm4}
\end{subtable}
\caption{Performance statistics by backtest scenario for each PM.}
\label{t-performance_stats_pm}
\end{table}

\end{document}